\documentclass[preprint]{elsarticle}

\usepackage{amsmath}
\usepackage{amssymb, amsfonts,amscd,verbatim,hyperref,cleveref,graphics}
\usepackage{xspace}
\usepackage{tikz}
\usetikzlibrary{matrix,calc,decorations.pathmorphing,shapes,arrows}
\newtheorem{thm}{Theorem}
\newtheorem{lemma}[thm]{Lemma}
\newtheorem{remark}[thm]{Remark}

\newcommand\numberthis{\addtocounter{equation}{1}\tag{\theequation}}

\def\bfm#1{\boldsymbol{#1}}
\def\NN{{\mathbb N}}

\def\RR{{\mathbb R}}

\DeclareSymbolFont{bbold}{U}{bbold}{m}{n}
\DeclareSymbolFontAlphabet{\mathbbold}{bbold}

\newcommand{\zs}

\begin{document}

\begin{frontmatter}
\title{Optimal parametric interpolants of circular arcs}

\author[FMF,IMFM]{Ale\v s Vavpeti\v c}
\ead{ales.vavpetic@fmf.uni-lj.si}
\address[FMF]{Faculty of Mathematics and Physics , University of Ljubljana, Jadranska 19, Ljubljana, Slovenia}
\address[IMFM]{Institute of Mathematics, Physics and Mechanics, Jadranska 19, Ljubljana, Slovenia}

\begin{abstract}
The aim of this paper is a construction of quartic parametric polynomial interpolants of a circular arc, where two boundary points of a circular arc are interpolated. For every unit circular arc of inner angle not greater than $\pi$ we find the best interpolant, where the optimality is measured by the simplified radial error.
\end{abstract}

\begin{keyword}
  geometric interpolation \sep circular arc \sep parametric polynomial \sep B\'ezier curve \sep optimal interpolation

\MSC 65D05 \sep 65D07 \sep 65D17
\end{keyword}

\end{frontmatter}

\section{Introduction}

Circular arcs are basic ingredients of several graphical and control systems, so their approximation by parametric polynomials is important in Computer Aided Geometric Design (CAGD), Computer Aided Design (CAD) and Computer Aided Manufacturing (CAM). Usually we construct parametric polynomial approximant of a circular arc by interpolation of some corresponding geometric quantities. This usually include interpolation of boundary points, corresponding tangent directions, signed curvatures... The results are so called geometric parametric polynomial interpolants (${\mathcal G}^n$ interpolants), which can be combined to form geometrically smooth spline curves.
One of the standard measures in this case is the radial distance $d_r$, measuring the distance of the point on the parametric polynomial to the corresponding point on the circular arc in the radial direction. Under some assumptions the metric $d_r$ is equivalent to the Hausdorff metric (\cite{Ahn-Kim-arcs-1997} and \cite{Jaklic-Kozak-18}). Hence to find the best interpolant of the unit circular arc $\bfm{c}$ with respect to the Hausdorff metric, we have to find an interpolant $\bfm{p}=(x,y)^T$ which minimizes the value $d_r(\bfm{c},\bfm{p})=\max_t\big|\|\bfm{p}(t)\|-1\big|=\max_t\big|\sqrt{x^2(t)+y^2(t)}-1\big|$.
In the very first paper in which the optimality was proved \cite{Morken-91}, M{\o}rken observed that the distance $d_r$ is rather cumbersome to work with, so he suggested that instead of the radial distance, we should use the simplified radial distance $d_{sr}$ defined by $d_{sr}(\bfm{c},\bfm{p})=\max_t\big|x^2(t)+y^2(t)-1\big|$. The involved function $x^2+y^2-1$ is a polynomial, which significantly simplifies an analysis of the optimality of the best interpolant.

There are many papers where different types of geometric approximations are considered. But only a few of them are dealing with the optimality of the solution. M{\o}rken considered the parabolic ${\mathcal G}^{0}$ interpolation of a circular arc \cite{Morken-91}. Hur and Kim analyzed the cubic ${\mathcal G}^{1}$ and the quartic ${\mathcal G}^{2}$ cases \cite{Hur-Kim-11}. In this three cases there is only one free parameter involved. Two parametric cases are considered by Vavpeti\v{c} and \v{Z}agar in \cite{Vavpetic-Zagar-19} where the optimal solutions for the cubic ${\mathcal G}^{0}$ and the quartic ${\mathcal G}^{1}$ interpolants were found. So far there are no results on the optimal solution of the quintic or higher degree interpolants.
\begin{table}[h]
\begin{center}
\begin{tabular}{cccc} \hline\hline
 order & 2 & 3 & 4 \\ \hline
${\mathcal G}^{0}$ & M{\o}rken (1991) & Vavpeti\v{c}, \v{Z}agar (2019) & this paper \\
${\mathcal G}^{1}$ & Knez, \v{Z}agar (2018) & Hur, Kim (2011) & Vavpeti\v{c}, \v{Z}agar (2019)  \\
${\mathcal G}^{2}$ & -- & Knez, \v{Z}agar (2018) & Hur, Kim (2011)  \\
${\mathcal G}^{3}$ & -- & -- & Knez, \v{Z}agar (2018) \\ \hline
\end{tabular}
\end{center}
\caption{The list of all results where the optimality of the solution was proved.}
\end{table}
In all cases the optimality is measured by simplified radial distance. The only paper where the optimality of the best interpolant is proved according to the real radial distance is \cite{Vavpetic-Zagar-preprint}.
Let us also mention the ${\mathcal G}^{n}$ interpolation of order $n+1$. For every circular arc and for every $n\in\NN$ there are only finitely many ${\mathcal G}^{n}$ interpolants of order $n+1$, and for many circular arcs there is only one such interpolant. Hence there is nothing much to optimize, the only question is the existence of the optimal interpolant, which was positively answered by Knez and \v{Z}agar \cite{Knez-Zagar-18}. In this paper we consider the only remaining case of order less then 5, i.e., the quartic ${\mathcal G}^{0}$ interpolation.

The paper is organised as follows. In Section~\ref{prelim} we review basic definitions and describe the idea of the construction of the best ${\mathcal G}^{0}$ interpolant of order $n$. In Sections~\ref{parabolic} and \ref{cubic} we use our method to confirm known results about the parabolic and the cubic interpolations of a circular arc. The main part of the paper is Section~\ref{quartic}, where we construct the best quartic interpolant of a circular arc and we prove its optimality according to the simplified radial distance. In Section~\ref{conc} we give some concluding remarks and suggestions for possible future research.

\section{Preliminaries}\label{prelim}

Let $0<\varphi\le\tfrac\pi 2$ and let $\bfm{c}\colon[-\varphi,\varphi]\to\RR^2$, $\bfm{c}(t)=(\cos t,\sin t)^T$, be the standard nonpolynomial parametrization of a unit circular arc. We'd like to find the best approximation of $\bfm{c}$ by polynomial curve $\bfm{p}\colon [-1,1]\to\RR^2$ of degree $n\in\NN$ for which $\bfm{p}(\pm 1)=(\cos \varphi,\pm\sin \varphi)^T$.
It is convenient to write $\bfm{p}=(x,y)^T$, where $x$ and $y$ are polynomials of degree at most $n$.
We shall choose the Bernstein-B\'ezier representation of $\bfm{p}$, i.e.,
\begin{equation}\label{p_Bern_Bez_form}
  \bfm{p}(t)=\sum_{j=0}^n B_j^n(t)\,\bfm{b}_j,
\end{equation}
where $B_j^n$, $j=0,1,\dots,n$, are (reparameterized) Bernstein polynomials over $[-1,1]$, given as
\begin{equation*}
  B_j^n(t)=\binom{n}{j}\left(\frac{1+t}{2}\right)^{j}\left(\frac{1-t}{2}\right)^{n-j},
\end{equation*}
and $\bfm{b}_j\in\RR^2$, $j=0,1,\dots,n$, are the control points. Since we consider ${\mathcal G}^0$ interpolation, we have $\bfm{b}_0=(\cos \varphi,-\sin \varphi)^T$ and $\bfm{b}_n=(\cos \varphi,\sin \varphi)^T$. The circular arc is symmetric over $x$ axis, therefore the best interpolant possesses the same symmetry, i.e., $\bfm{b}_{n-j}=r(\bfm{b}_j)$ for all $j=0,1,\dots,n$, where $r\colon\RR^2\to\RR^2$ is the reflection over $x$ axis. Therefore all possible sets of control points of desired interpolants can be described by $n-1$ parameters.

\begin{figure}[h]
\begin{center}
\includegraphics[width=0.33\textwidth]{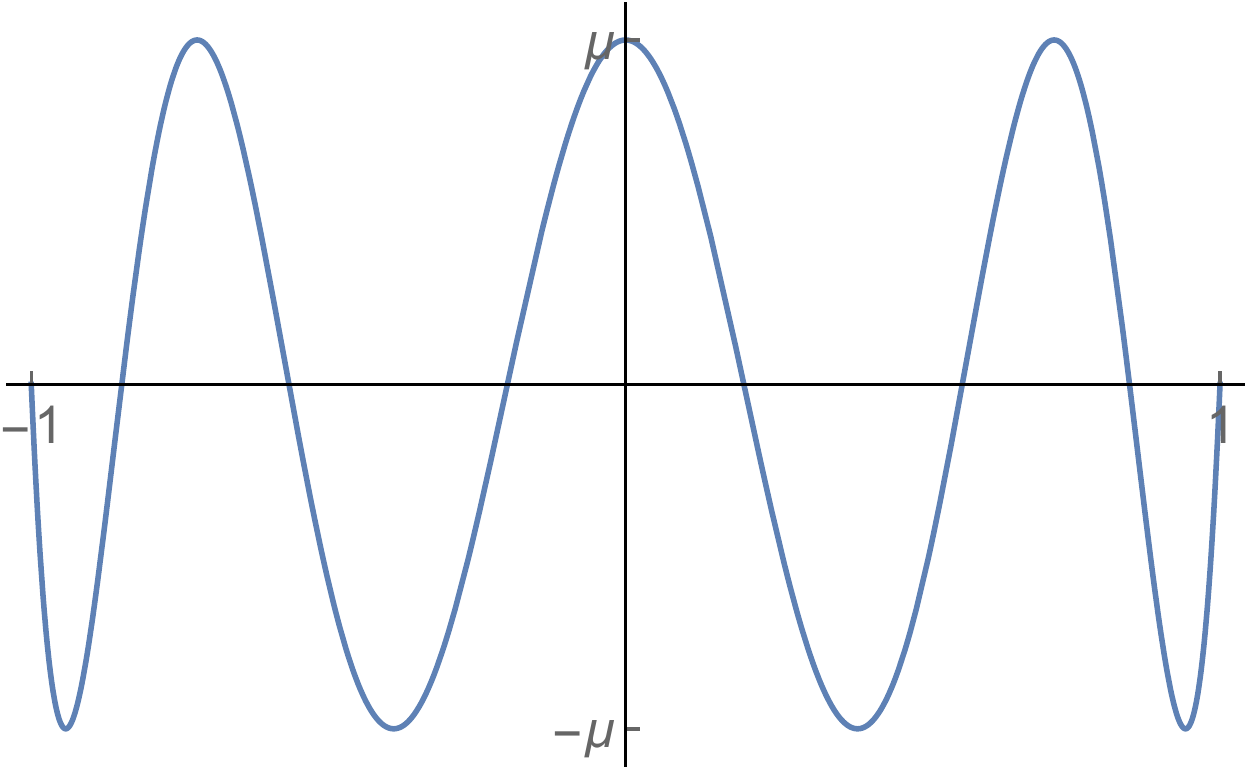}
\end{center}
\caption{The graph of the error function of the best quartic ${\mathcal G}^0$ interpolant.}\label{graph_error}
\end{figure}

The simplified signed radial error function $\psi$ will be defined as
\begin{align*}
  \psi(t)&=x^2(t)+y^2(t)-1=\|\bfm{p}(t)\|_2^2-1,\quad t\in[-1,1],
\end{align*}
where $\|\cdot\|_2$ is the Euclidean norm. In the case of ${\mathcal G}^0$ interpolant, the corresponding simplified error function has zeros at $\pm 1$.
Eisele showed \cite{Eisele-94} that the best ${\mathcal G}^k$ interpolant of order $n$
(if it exists) is an alternant with $2(n-k-1)+1$ extreme points, i.e., the corresponding signed simplified error function has $2(n-k-1)+1$ local extrema of the same absolute value and sequential ones have different sign (see Figure~\ref{graph_error}).
In the case of ${\mathcal G}^0$ interpolants which is considered in this paper, the simplified error function of the best interpolant of order $n$ is of the form $\chi(t)=\mu T_{2n}(\zeta_{n} t)$, where $T_{2n}$ is the Chebyshev polynomial of order $2n$, $\zeta_{n}=\cos\tfrac{\pi}{2n}$ is its largest zero and $\mu\in\RR\setminus\{0\}$ is a multiplicative constant.

The candidates for the best interpolants are those which have the coresponding simplified error function of the form $\chi$. Hence, to find all candidates for the best interpolant of a circular arc of order $n$ we have to solve the system of equations $\psi(\zeta_i)=0$, where $\zeta_i=(\cos\tfrac{\pi}{2n})^{-1}\cos(\tfrac{2i+1}{2n}\pi)$, $i=1,\ldots,n-1$, are zeros of the function $\chi$ on the interval $(0,1)$. This is a nonlinear system of $n-1$ polynomial equations for $n-1$ parameters (which describe control points of an interpolant). The system is solved using Gr{\"o}bner basis. Namely, we find a Gr{\"o}bner basis for the ideal $\langle \psi(\zeta_i)\mid i=1,\ldots,n-1\rangle$ for some monomial order, such that one element of the basis is a polynomial $f$ in only one variable (of order $n$). All interpolants obtained as a solution of the above system of equations induce error functions of the same shape (see Figure \ref{graph_error}), which differ only by the multiplicative constant $\mu$. So we have to find the zero of $f$, such that the absolute value of the corresponding multiplicative constant $\mu$ is the smallest possible, or equivalently, such that the absolute value of the leading coefficient (or any other coefficient) of the simplified error function is the smallest possible.

In what follows many mathematical expressions are considered as polynomial in $\cos\varphi$, so it makes sense to define new variables $c:=\cos\varphi$ and $s:=\sin\varphi$.

\section{The parabolic ${\mathcal G}^{0}$ case}\label{parabolic}

The first paper, where polynomial interpolants of circular arcs are considered in Bernstein-B\'ezier form, is \cite{Morken-91}. In this paper M{\o}rken found the best parabolic ${\mathcal G}^{0}$ interpolant for arbitrary circular arc with inner angle not greater then $\pi$. Using Eisele's theorem \cite{Eisele-94} it is now very easy to find such interpolant. The control points of an arbitrary parabolic ${\mathcal G}^0$ interpolant are $\bfm{b}_0=(c,-s)^T$, $\bfm{b}_1=(d,0)^T$, $\bfm{b}_2=(c,s)^T$, where $d\in\RR$ is an unknown parameter, and the corresponding simplified error function is
$$
\psi(t)=\frac {1} 4 \left(t^2-1\right)\left((d-c)^2t^2-(d+c)^2+4\right).
$$
Because the simplified error function of the best interpolant must be of the form $\chi(t)=\mu T_{4}(\zeta_{2} t)=c(t^2-1)(t^2-(3-2\sqrt{2}))$, we get the equation
$$
\frac{(d+c)^2-4}{(d-c)^{2}}=(3-2\sqrt{2})
$$
for the parameter $d$, which gives the best interpolant. The only positive solution is
$$
d=\left(1+\sqrt{2}\right) \left(\sqrt{2 \left(\sqrt{2}-1\right)+\left(3-2 \sqrt{2}\right) c^2}-\left(2-\sqrt{2}\right) c\right).
$$
\vskip-8mm
\begin{table}[h]
\begin{center}
$$
\begin{array}{ccccc} \hline\hline
\varphi && d &&  \text{radial error}\\ \hline
\frac\pi2 && 2.19737 && 2.07107\times 10^{-1} \\
\frac\pi3 && 1.54643 && 4.69687\times 10^{-2} \\
\frac\pi4 && 1.30834 && 1.55050\times 10^{-2} \\
\frac\pi6 && 1.13712 && 3.15242\times 10^{-3} \\
\frac\pi8 && 1.07713 && 1.00735\times 10^{-3} \\
\frac\pi{12}&& 1.03427 && 2.00378\times 10^{-4} \\ \hline
\end{array}
$$
\end{center}
\caption{A table of the simplified radial errors of the best parabolic ${\mathcal G}^0$ geometric interpolants of a circular arc given by the inner angle of $2\varphi$. }
\end{table}

\section{The cubic ${\mathcal G}^{0}$ case}\label{cubic}

This case was considered in \cite{Vavpetic-Zagar-19}. A little different approach is presented here, which was described in Section~\ref{prelim} and which can be generalised to the quartic ${\mathcal G}^{0}$ case. The control points for an arbitrary cubic ${\mathcal G}^0$ interpolant are $\bfm{b}_0=(c,-s)^T$, $\bfm{b}_1=(\xi,-\eta)^T$, $\bfm{b}_2=(\xi,\eta)^T$ $\bfm{b}_3=(c,s)^T$, and the corresponding signed simplified error function is
\begin{align*}
\psi(t)=&\frac{1}{16} (t^2-1) \big((3 \eta -s)^2 t^4+(16 s^2-9 (\eta +s)^2+9 (\xi -c)^2) t^2
+(16-(3 \xi +c)^2)\big).
\end{align*}
The simplified error function of the best interpolant is of the form $\chi(t)=\mu T_{6}(\zeta_{3} t)=c(t^2-1)(t^2-u^2)(t^2-v^2)$, where $u=\sqrt{3}-1$ and $v=2-\sqrt{3}$. To get parameters $(\xi,\eta)$ of the best interpolant, we have to solve the system of equations: $\psi(u)=\psi(v)=0$.
From $(1-v^2)^2v^2\psi(u)-(1-u^2)^2u^2\psi(v)=0$, we get
\begin{align}\label{cubic1}
\eta=\frac1{s}\left(\frac{2+\sqrt{3}}{8} \left(3 \xi +c\right)^2-\xi  c-3-2\sqrt{3}\right).
\end{align}
To keep the symmetry we use \eqref{cubic1} in the equality $\psi(u)-\psi(v)=0$. We get one solution, $\xi=c$, and the remaining ones satisfy the equality
%\begin{align}\label{equationcubic}
%f(\xi):=&243 \xi ^3-27 c \left(11-16 \sqrt{3}\right) \xi ^2-3 \left(32\left(1+2\sqrt{3}\right)-3 \left(81-32\sqrt{3}\right) c^2\right) \xi \nonumber\\
%&-32\left(13+2\sqrt{3}\right) c-\left(163-112 \sqrt{3}\right) c^3=0.
%\end{align}
\begin{align}\label{equationcubic}
f(\xi):=243 \xi ^3-27 c \left(11-16 \sqrt{3}\right) \xi ^2-3 \left(32\left(1+2\sqrt{3}\right)-3 \left(81-32\sqrt{3}\right) c^2\right) \xi
-32\left(13+2\sqrt{3}\right) c-\left(163-112 \sqrt{3}\right) c^3=0.
\end{align}
Because $f(0)=-(579-48\sqrt 3)c-(112\sqrt 3-163)cs^2<0$, $f(-\tfrac 1 9(8\sqrt 3-1)c)=\tfrac{256}3cs^2>0$,
and the leading coefficient of the cubic polynomial $f$ is positive, $f$ has only one positive solution. Since $f(c)=-256(2+\sqrt 3)c s^2<0$ and
$f(\tfrac 1 3(4+c))=8(56-32\sqrt 3+(80\sqrt 3-68)c+102c^2+(7+8\sqrt 3)c^3)>0$, the positive solution is on the interval $(c,\tfrac 1 3(4+c))$.
We have two candidates for the best interpolant, induced by $\xi_1=c$ and $\xi_2\in(c,\tfrac 1 3(4+c))$. Both candidates induce an error function of the form $\psi(t)=\mu T_{6}(\zeta_{3} t)$. Therefore we want to minimize the absolute value of the constant $\mu$, which is equivalent to minimizing the absolute value of the constant coefficient of the simplified error polynomial $\psi$. The constant coefficient of the polynomial $\psi$ is $-\tfrac 1{16}(16-(3\xi+c)^2)$, which is negative and increasing function of the variable $\xi$ on the interval $[c,\tfrac 1 3(4+c)]$, hence the best interpolant is induced by the parameter $\xi_2$. We proved the following theorem.

\begin{thm}
For every $\varphi\in (0,\tfrac\pi 2]$ there exists the unique pair of parameters $(\xi,\eta)$ which induces the best polynomial interpolant of a given circular arc. The parameter $\xi$ is the only zero of the function \eqref{equationcubic} on the interval $(\cos\varphi,\tfrac 1 3(4+\cos\varphi))$ and $\eta$ is defined by \eqref{cubic1}.
\end{thm}

\begin{remark}
A purely geometric observation reveals that the solution $\xi=c$ is not an admissible one. The first coordinate of the cubic interpolant \eqref{p_Bern_Bez_form} is $x(t)=\tfrac1 4(3\xi+c-3t^2(\xi-c))$. For $\xi=c$ the function $x(t)=c$ is a constant, therefore the interpolant is a line segment which is not a desired solution. Hence also in the cubic case for every $\varphi$ we get only one admissible candidate for the best interpolant of the circular arc with the inner angle $2\varphi$.
\end{remark}

\vspace{-8mm}
\begin{table}[h]
\begin{center}
$$
\begin{array}{ccccccc} \hline\hline
\varphi && \xi && \eta && \text{radial error}\\ \hline
\frac\pi2 && 1.32800 && 0.94046 && 7.97742\times 10^{-3} \\
\frac\pi3 && 1.16617 && 0.47494 && 7.50902\times 10^{-4} \\
\frac\pi4 && 1.09754 && 0.31523 && 1.36878\times 10^{-4} \\
\frac\pi6 && 1.04465 && 0.19043 && 1.22221\times 10^{-5} \\
\frac\pi8 && 1.02537 && 0.13762 && 2.18815\times 10^{-6} \\
\frac\pi{12}&& 1.01136 && 0.08926 && 1.92912\times 10^{-7} \\ \hline
\end{array}
$$
\end{center}
\caption{A table of the simplified radial errors of the best cubic ${\mathcal G}^0$ geometric interpolants of a circular arc given by the inner angle of $2\varphi$.}
\end{table}

\section{The quartic ${\mathcal G}^{0}$ case}\label{quartic}

This is a tree-parametric problem.
The control points are $\bfm{b}_0=(c,-s)^T$, $\bfm{b}_1=(\alpha,\beta)^T$, $\bfm{b}_2=(\gamma,0)^T$, $\bfm{b}_3=(\alpha,-\beta)^T$, $\bfm{b}_4=(c,s)^T$, and the corresponding signed simplified error function is
\begin{align*}
\psi(t)=-1+\frac{1}{64} \left(4 \left(1-t^4\right) \alpha +3 \left(1-t^2\right)^2 \gamma +\left(1+6 t^2+t^4\right) c \right)^2
+\frac{1}{4} t^2\left(2 \left(1-t^2\right) \beta +\left(1+t^2\right) s\right)^2.
\end{align*}
By Eisele's theorem the simplified error function of the best interpolant is of the form $\chi(t)=\mu T_8\left(\zeta_4 t\right)$, with three zeros $u_1=\sqrt{2 \left(2+\sqrt{2}\right)}-1-\sqrt{2}$, $u_2=\sqrt{2+\sqrt{2}}-1$ and $u_3=1+\sqrt{2}-\sqrt{2+\sqrt{2}}$ on the interval $(0,1)$.
We have to solve the system of equations $\psi(u_j)=0$, $j=1,2,3$, and find out
which solution $(\alpha,\beta,\gamma)$ induces the best interpolation of the circular arc. In what follows it is useful to define $\sigma_1=(1-u_1^2)+(1-u_1^2)+(1-u_1^2)\approx 1.92$, $\sigma_2=(1-u_1^2)(1-u_2^2)+(1-u_1^2)(1-u_3^2)+(1-u_2^2)(1-u_3^2)\approx 1.11$, and $\sigma_3=(1-u_1^2)(1-u_1^2)(1-u_1^2)\approx 0.18$.

We form the linear combination of equations of the system so that we eliminate the variable $\beta$ and get
\begin{align*}
0&=64\left(\frac{u_{2}^2u_{3}^2\psi(u_1)}{(u_{2}^2-u_1^2)(u_1^2-u_{3}^2)(1-u_1^2)}+
\frac{u_{1}^2u_{3}^2\psi(u_2)}{(u_{3}^2-u_2^2)(u_2^2-u_{1}^2)(1-u_2^2)}+
\frac{u_{1}^2u_{2}^2\psi(u_3)}{(u_{1}^2-u_3^2)(u_3^2-u_{2}^2)(1-u_3^2)}\right)\\
&=64+u_1^2 u_2^2 u_3^2 (4 \alpha -3 \gamma -c)^2-(4 \alpha +3 \gamma +c)^2,
\numberthis\label{quartic1}
\end{align*}
%\begin{align*}
%0&=64\sum_{j=1}^3\frac{u_{j-1}^2u_{j+1}^2\psi(u_j)}{(u_{j+1}^2-u_j^2)(u_j^2-u_{j-1}^2)(1-u_j^2)}\\
%&=64+u_1^2 u_2^2 u_3^2 (4 \alpha -3 \gamma -c)^2-(4 \alpha +3 \gamma +c)^2,
%\numberthis\label{quartic1}
%\end{align*}
then set $x:=4 \alpha -3 \gamma -c$ and $y:=4 \alpha +3 \gamma +c$. Similarly we eliminate $\beta^2$ and $y^2$ and get
\begin{align*}
0&=\sum_{j=1}^3(u_{j-1}^2 \left(1-u_{j-1}^2\right)^2-u_{j+1}^2 \left(1-u_{j+1}^2\right)^2)\psi(u_j)
   -\frac{\sigma_2}{64}\left(u_1^2-u_2^2\right) \left(u_2^2-u_3^2\right)\left(u_3^2-u_1^2\right) \left(64+u_1^2 u_2^2 u_3^2 x^2-y^2\right)\\
&=2\left(u_1^2-u_2^2\right) \left(u_2^2-u_3^2\right)\left(u_3^2-u_1^2\right)(\sigma_2-\sigma_1+1)\left(-1+\frac{1}{128} \sigma_3 x^2+\frac{1}{8} (x+y) c+\beta s\right),
\end{align*}
where $u_0=u_3$ and $u_4=u_1$, therefore
\begin{align}\label{beta_quartic}
0=-1+\frac{1}{128} \sigma_3 x^2+\frac{1}{8} (x+y) c+\beta s.
\end{align}
Using the equalities \eqref{quartic1} and \eqref{beta_quartic} we get
\begin{align*}
0&=\frac{\psi(u_1)}{1-u_1^2}+\frac{\psi(u_2)}{1-u_2^2}+\frac{\psi(u_3)}{1-u_3^2}
+\frac1{64}\sigma_1(64+u_1^2u_2^2u_3^2x^2-y^2)
-\left(\sigma_1^2-3\sigma_1-2\sigma_2+6\right)\left(-1+\frac{1}{128} \sigma_3 x^2+\frac{1}{8}c (x+y)+\beta s\right)\\
&=\frac{1}{128} \left(\sigma_1^2-2 \sigma_2-\sigma_1\right) \left(96-4 x y-128 \beta ^2+
   16c (x-y)+ 32c^2- \left(4-2\sigma_1+\sigma_3\right)x^2\right),
\end{align*}
so
\begin{align*}
0&=96-4 x y-128 \beta ^2+
   16 c(x-y)+ 32c^2- \left(4-2\sigma_1+\sigma_3\right)x^2.
\end{align*}
We multiply the last equality by $s^2$, use the equality $\beta^2s^2=(-1+\frac{1}{128} \sigma_3 x^2+\frac{1}{8}c (x+y))^2$ obtained from \eqref{beta_quartic}, then use the equality $y^2=64+u_1^2u_2^2u_3^2x^2=64+(1-s_1+s_2-s_3)x^2$ obtained from \eqref{quartic1},
and get
%\begin{align*}
%0&=\left(-\frac{x}{32}+\frac{c}{8}+\frac{c^3}{8}-\frac{1}{512} \sigma_3c x^2 \right)y
%-\frac1 4\bigg(\left(\frac{1}{64} \sigma_3 x^2+\frac{1}{4} cx-s^2\right)^2\\
%&-\frac{1}{16} \left(1-\sigma_2+2 \sigma_3\right)  c^2x^2+\frac{1}{16} \left(2-\sigma_1\right) x^2-c (x-8 c)\bigg). \numberthis\label{y_quartic}
%\end{align*}
\begin{align}\label{y_quartic}
0=\left(-\frac{x}{32}+\frac{c}{8}+\frac{c^3}{8}-\frac{1}{512} \sigma_3c x^2 \right)y
-\frac1 4\bigg(\left(\frac{1}{64} \sigma_3 x^2+\frac{1}{4} cx-s^2\right)^2
-\frac{1}{16} \left(1-\sigma_2+2 \sigma_3\right)  c^2x^2+\frac{1}{16} \left(2-\sigma_1\right) x^2-c (x-8 c)\bigg).
\end{align}
By combining the equalities \eqref{y_quartic} and \eqref{quartic1}, we see that we have to investigate the zeros of the function
\begin{align*}
f(x)&:=\frac{1}{16}\biggl(\!\biggl(\frac{1}{64} \sigma_3 x^2+\frac{1}{4} c x-s^2\biggr)^2
-\frac{1}{16} \left(1-\sigma_2+2 \sigma_3\right) c^2 x^2
+\frac{1}{16} \left(2-\sigma_1\right) x^2- c(x-8 c) \!\biggr)^2\\
&-(64+u_1^2u_2^2u_3^2x^2)\left(\frac{x}{32}-\frac{c}{8}-\frac{c^3}{8}+\frac{1}{512} \sigma_3c x^2 \right)^2.\numberthis\label{function_quartic}
\end{align*}
\begin{figure}[h]
\begin{center}
\includegraphics[width=0.8\textwidth]{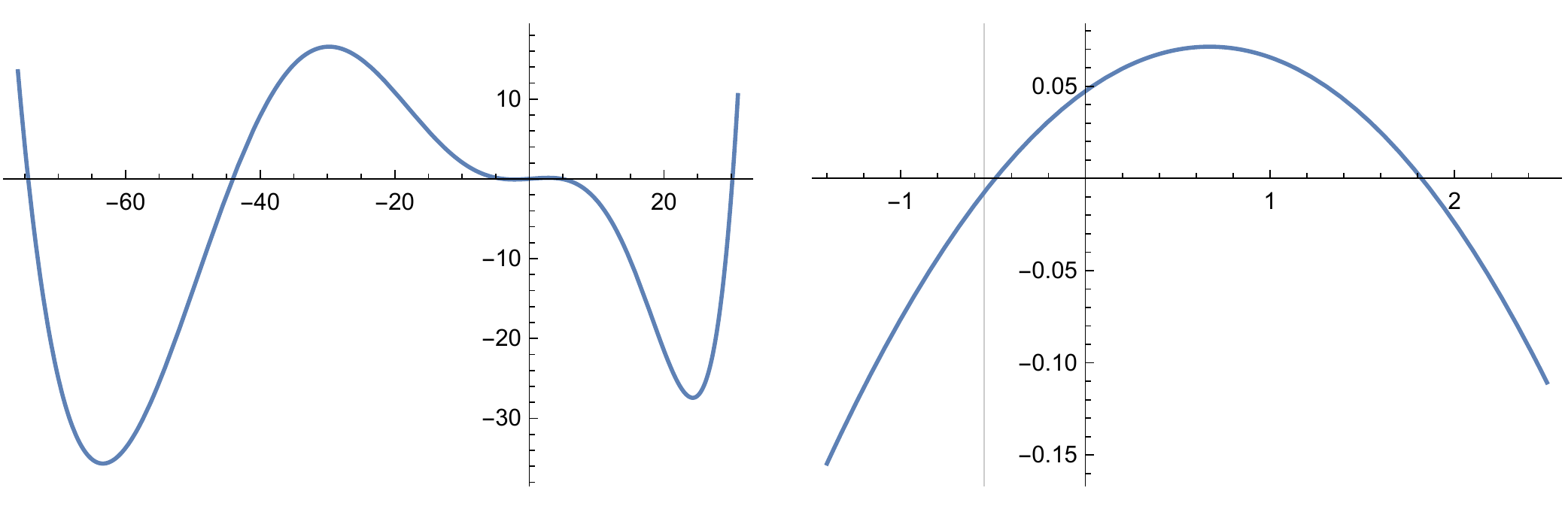}
\caption{The graph of the function $f$ for the angle $\varphi=\tfrac 1 4 \pi$ is on the left and the graph of the function $f$ for the angle $\varphi=\tfrac{5}{12} \pi$ on smaller interval around 0 where the desired zero of $f$ appears on the right. The vertical line on the right graph is drawn at $x=-(1-\cos\varphi)^2$.}
\end{center}
\end{figure}
Note that for every real zero $x$ of the function $f$ there are the unique (real) number $y$ obtained from \eqref{y_quartic} and the unique (real) number $\beta$ obtained from \eqref{beta_quartic}, such that the triple $(\alpha=\tfrac 1 8(x+y),\beta,\gamma=\tfrac 1 6(y-x-2\cos\varphi))$ induces an interpolant with the simplified error function $\psi$ of the form $\chi$.
The amplitude of the simplified error function $\psi$ is the smallest possible if the leading coefficient of $\psi$ is the smallest possible. The leading coefficient of $\psi$ is $\tfrac{x^2}{64}$, hence we must find the zero of $f$ with the smallest absolute value.

\begin{remark}\label{positive}
Quite often we will use the following argument:
Let $p(x)=c_0+c_1x+\ldots+c_nx^n$ be a nonzero polynomial such that the sum $d_j:=c_0+\ldots+c_j$ is nonnegative for all $j=0,\ldots,n$. Then we can write $p(x)=(1-x)(d_0+d_1x+\ldots+d_{n-1}x^{n-1})+d_nx^n$, therefore $p(x)> 0$ for all $x\in(0,1)$. Similarly, if $d_j\le 0$ for all $j=0,\ldots,n$, then $p(x)< 0$ for all $x\in(0,1)$.
\end{remark}

\begin{lemma}
The function $f$ has a zero on the interval $[-(1-c)^2,0]$.
\end{lemma}

\begin{pf}
Note that $f(0)=\tfrac{s^8}{16}>0$. We can write
$f(-(1-c)^2)=-\left(f_1^2(c)+(1-c)^4f_2(c)f_3(c)\right)$, where
\begin{align*}
f_1(c)&=\frac 1 {512}u_1u_2u_3(1-c)^2\Big(
\sigma_3(1-c)^4c-16(1+2c+c^2+4c^3)\Big),\\
f_2(c)&=
-\frac{1}{32}+\frac{\sigma_1}{64}+\frac{\sigma_3}{128}-\frac{\sigma_3^2}{16384}+
\left(\frac{1}{8}-\frac{7\sigma_3}{512}+\frac{\sigma_3^2}{4096}\right) c
+\left(-\frac{\sigma_2}{64}+\frac{5\sigma_3}{256}-\frac{3 \sigma_3^2}{8192}\right)c^2
+\left(\frac{\sigma_3}{512}+\frac{\sigma_3^2}{4096}\right) c^3
-\frac{\sigma_3^2}{16384} c^4,\\
f_3(c)&=
\frac{17}{32}-\frac{\sigma_1}{64}-\frac{\sigma_3}{128}+\frac{\sigma_3^2}{16384}+
\left(\frac{3}{4}+\frac{\sigma_1}{16}+\frac{7\sigma_3}{512}-\frac{\sigma_3^2}{2048}\right) c
+\left(\frac{19}{16}-\frac{3\sigma_1}{32}+\frac{\sigma_2}{64}+\frac{\sigma_3}{256}+\frac{7\sigma_3^2}{4096}\right)
c^2\\
&+\left(\frac{9}{8}+\frac{\sigma_1}{16}-\frac{\sigma_2}{16}+\frac{\sigma_3}{512}-\frac{7 \sigma_3^2}{2048}\right)c^3
+\left(\frac{17}{32}-\frac{\sigma_1}{64}+\frac{3\sigma_2}{32}-\frac{3 \sigma_3}{64}+\frac{35 \sigma_3^2}{8192}\right)c^4
+\left(-\frac{1}{8}-\frac{\sigma_2}{16}+\frac{25\sigma_3}{512}-\frac{7 \sigma_3^2}{2048}\right) c^5\\
&+\left(\frac{\sigma_2}{64}-\frac{3 \sigma_3}{256}+\frac{7\sigma_3^2}{4096}\right)c^6
+\left(-\frac{\sigma_3}{512}-\frac{\sigma_3^2}{2048}\right) c^7
+\frac{\sigma_3^2}{16384} c^8.
\end{align*}
By Remark~\ref{positive} it is easy to see that $f_2(c),f_3(c)>0$ for all $c\in(0,1]$, hence $f(-(1-c)^2)< 0$, therefore $f$ has a zero on the interval $[-(1-c)^2,0]$.
\qed
\end{pf}

\begin{lemma}
Let $\zeta_-$ be the largest negative zero of the function $f$. If $\varphi<\tfrac \pi 2$, then for every positive zero $\zeta_+$ of $f$, we have $\zeta_+>|\zeta_-|$. If $\varphi=\tfrac\pi 2$, then $f$ has a positive zero and the smallest positive zero $\zeta_+$ of $f$ satisfies $\zeta_+=|\zeta_-|$.
\end{lemma}

\begin{pf}
If $\varphi=\tfrac\pi 2$, then $f$ is an even function, hence $|\zeta_-|=\zeta_+$.

Let $\varphi<\tfrac\pi 2$.
By the previous lemma we know that $f$ has a zero on the interval $[-(1-c)^2,0]$ and $f(0)>0$. So it is enough to show that $f(x)> f(-x)$ for all $x\in(0,(1-c)^2]$. Let us define $g(x)=f(x)-f(-x)$. Since $g(0)=0$ it is enough to prove that the function $g$ is convex on $[0,(1-c)^2]$ and $g((1-c)^2)> 0$.

The second derivative of $g$ is
\begin{align*}
g''(x)&=\frac{c}{256}x\left(g_1(c)+ g_2(c)x^2+\frac{21 \sigma_3^3}{4096} x^4\right),
\end{align*}
where
\begin{align*}
g_1(c)&=3 \left(\left(-16+4 \sigma_1+8 \sigma_2-9 \sigma_3\right)+\left(16-12 \sigma_1-4 \sigma_2+14 \sigma_3\right) c^2+\left(4 \sigma_2-5 \sigma_3\right) c^4\right),\\
g_2(c)&=\frac{5}{32} \sigma_3 \left(4 \sigma_1-8 \sigma_2+3 \sigma_3+\left(4 \sigma_2-5 \sigma_3\right) c^2\right).
\end{align*}
By Remark~\ref{positive}, $g_1(c)<0$, $g_1(c)+g_2(c)<0$ and $g_1(c)+g_2(c)+\tfrac{21}{4096}\sigma_3^3<0$ for all $c\in[0,1]$. Then, again by Remark~\ref{positive}, we get $g''(x)< 0$ for all $x\in(0,1)$ and all $c\in[0,1]$, hence $g$ is a convex function on $[0,1]$ for all $c\in [0,1]$.

Since
\begin{align*}
g((1-c)^2)=c(1-c)^4(\frac 1 3(1+2c)^2+\frac 1{64}(1-c)^2g_3(c)),
\end{align*}
where
\begin{align*}
g_3(c)&=\frac{1}{3} \left(176-3 \sigma_1+3 \sigma_2\right)-\left(52-3 \sigma_1+\sigma_2+\sigma_3\right)(1-c)+\frac{1}{2} \left(20-3 \sigma_1+5 \sigma_2-4 \sigma_3\right) (1-c)^2\\
&-\frac{1}{2} \left(4 \sigma_2-5 \sigma_3\right) (1-c)^3+\frac{1}{256} \left(2 \sigma_2 \left(64-\sigma_3\right)-\left(160-2 \sigma_1+\sigma_3\right) \sigma_3\right) (1-c)^4\\
&-\frac{1}{256} \sigma_3 \left(4 \sigma_2-5 \sigma_3\right) (1-c)^5+\frac{1}{512}\left(4 \sigma_2-5 \sigma_3\right) \sigma_3(1-c)^6+\frac{1}{32768}\sigma_3^3 (1-c)^8
\end{align*}
it is enough to prove that $g_3(c)>0$ for all $c\in [0,1]$. The function $g_3$ is a polynomial in variable $(1-c)$ and by Remark~\ref{positive}, $g_3(c)> 0$ for all $c\in[0,1]$.
\qed
\end{pf}

By the previous lemma it seems that for $\varphi=\tfrac\pi 2$, there are two candidates for the best interpolation, but the next lemma shows that one is not admissible.

\begin{lemma}
Let $\varphi=\tfrac\pi 2$ and let $\zeta_+$ be the smallest positive zero of the function $f$. The interpolant induced by $\zeta_+$ is not admissible.
\end{lemma}

\begin{pf}
By \eqref{y_quartic} we get
\begin{align*}
y&=-\frac{32}{\zeta_+}\left(\frac 1{64}\zeta_+^2(2-\sigma_1)+\left(\frac12-\frac{\zeta_+^2\sigma_3}{128}\right)^2\right)
<-\frac{32}{\zeta_+}\left(\frac12-\frac{1}{128}\right)^2<-\frac{63^2}{16\cdot 32}<-1.
\end{align*}
Then $\alpha=\tfrac 1 8(x+y)=\tfrac1 8(\zeta_++y)<0$ and $\gamma=\tfrac 1 6(y-x-2c)\le\tfrac1 6(y-\zeta_+)<0$,  hence by the convex hull property, the whole B\'ezier polygon lies left from $y$-axis, therefore the interpolant is not a desired one.
\qed
\end{pf}

We proved that the best interpolant is induced by the largest negative zero of the function $f$ on $[-(1-c)^2,0]$. For numerical computations it would be desired that $f$ has only one zero on that interval.

\begin{lemma}
There is exactly one zero of $f$ on the interval $[-(1-c)^2,0]$.
\end{lemma}

\begin{pf}
It is enough to show that $f'(x)>0$ for all $x\in(-(1-c)^2,0)$ which is equivalent to $g'(t)<0$ for all $t\in(0,1)$, where $g(t)=f(-t(1-c)^2)$. We can write $g'(t)=-\left((1-c)^4 g_1(t)+\tfrac1{8192}(1-c)^8\left(g_2(t)+g_3(t)+g_4(t)\right)\right)$ where
\begin{align*}
g_1(t)&=\frac{1}{16} (1+c)^2 \left(5+c^2\right)c
+\frac{1}{128} \left(12-56 c^2+28 c^4+2\sigma_1 \left(1-c^2\right) \left(1+5 c^2+2 c^4\right)+2 \sigma_2
   \left(1-c^2\right)^2c^2+ \sigma_3\left(1-c^2\right)^3\right) t\\
&+\frac{3}{1024} \left(4 \sigma_1\left(1-3 c^2\right)+4 \sigma_2\left(2-c^2+c^4\right) + \left(1-c^2\right)\left(-16+\left(-9+5 c^2\right)\sigma_3\right)\right) (1-c)^2 c t^2,\\
g_2(t)&=\left(\frac{5}{2}\sigma_3-\sigma_2 c\right)^2 t^3,\\
g_3(t)&=\Bigl(\sigma_2 \left(32\left(1-c^2\right)+
     \sigma_3 \left(5-8 c+8 c^3\right)c\right)-\frac{1}{4} \sigma_3 \left(64 \left(1-5 c^2+2 c^4\right)+\sigma_3\left(37-24 c^2+12 c^4\right)\right)+8 \sigma_1 \left(2\sigma_2 c^2
   - \sigma_3\left(1+c^2-2 c^4\right)\right)\\
&-8 \sigma_1^2-\left(c^2+8 c^4\right) \sigma_2^2\Bigr)t^3
-\frac{5}{8} \sigma_3 \left(8 \sigma_2+5 c^2 \sigma_3\right)(1-c)^2 c t^4 -\frac{3}{64} \sigma_3^2 \left(4+4 c^2 \sigma_1+c^2
   \sigma_3\right)(1-c)^4 t^5-\frac1{4096} \sigma_3^4(1-c)^8t^7,\\
g_4(t)&=\frac{5}{8} \sigma_3 \left(4 \sigma_1+4 c^2 \sigma_2+3 \sigma_3\right)(1-c)^2 c t^4+\frac{3}{64} \sigma_3^2 \left(2 \sigma_1+2 c^2 \sigma_2+\sigma_3\right) (1-c)^4 t^5
+\frac{7}{512} \sigma_3^3(1-c)^6 c  t^6.
\end{align*}
It is easy to see that $g_2(t)> 0$ and $g_4(t)\ge 0$ for all $t\in (0,1)$ and $c\in [0,1)$. By Remark~\ref{positive}, also $g_1(t)\ge 0$ for all $t\in [0,1]$ and $c\in [0,1)$. The coefficients of $t^4$, $t^5$ and $t^7$ in the polynomial $g_3$ are obviously negative for all $c\in[0,1)$. Hence if we replace all $t^k$ in $g_3$ by $t^3$, we decrease the value of $g_3(t)$ and we get the expression of the form $t^3 h(c)$, where $h$ is a polynomial. By Remark~\ref{positive} we get $h(c)>0$ for all $c\in [0,1)$, hence $g_2(t)>0$ for all $t\in(0,1)$ and $c\in [0,1)$. This implies that $f$ is increasing on the interval $[-(1-c)^2,0]$.
\qed
\end{pf}

We proved the following theorem.

\begin{thm}
For every $\varphi\in (0,\tfrac\pi 2]$ there exists the unique triple of parameters $(\alpha,\beta,\gamma)$, which induces the best polynomial interpolant of a circular arc of the inner angle $2\varphi$, where $\alpha=\tfrac1 8 (x+y)$, $\gamma=\tfrac1 6(y-x-2\cos\varphi)$, $x$ is the only zero of the function \eqref{function_quartic} on the interval $[-(1-\cos\varphi)^2,0]$, $y$ is the solution of the equation \eqref{y_quartic}, and $\beta$ is the solution of the equation \eqref{beta_quartic}.
\end{thm}

\begin{figure}[h]\label{quartic_table}
\begin{center}
$$
\begin{array}{ccccc} \hline\hline
\varphi & \alpha & \beta & \gamma & \text{radial error}\\ \hline
\frac\pi2 & 0.87518 & 0.99857 & 1.49995 & 1.42325\times 10^{-4} \\
\frac\pi3 & 0.97471 & 0.59188 & 1.20039 & 5.83570\times 10^{-6} \\
\frac\pi4 & 0.99193 & 0.42228 & 1.10839 & 5.94378\times 10^{-7} \\
\frac\pi6 & 0.99840 & 0.27073 & 1.04680 & 2.34778\times 10^{-8} \\
\frac\pi8 & 0.99949 & 0.20014 & 1.02605 & 2.36051\times 10^{-9} \\
\frac\pi{12}& 0.99990 & 0.13203 & 1.01149 & 9.23852\times 10^{-11} \\ \hline
\end{array}
$$
\end{center}
\caption{A table of the radial errors of the best quartic ${\mathcal G}^0$ geometric interpolants of an circular arc given by inner angle of $2\varphi$.}
\end{figure}

\begin{figure}[h]
\begin{center}
\includegraphics[width=0.8\textwidth]{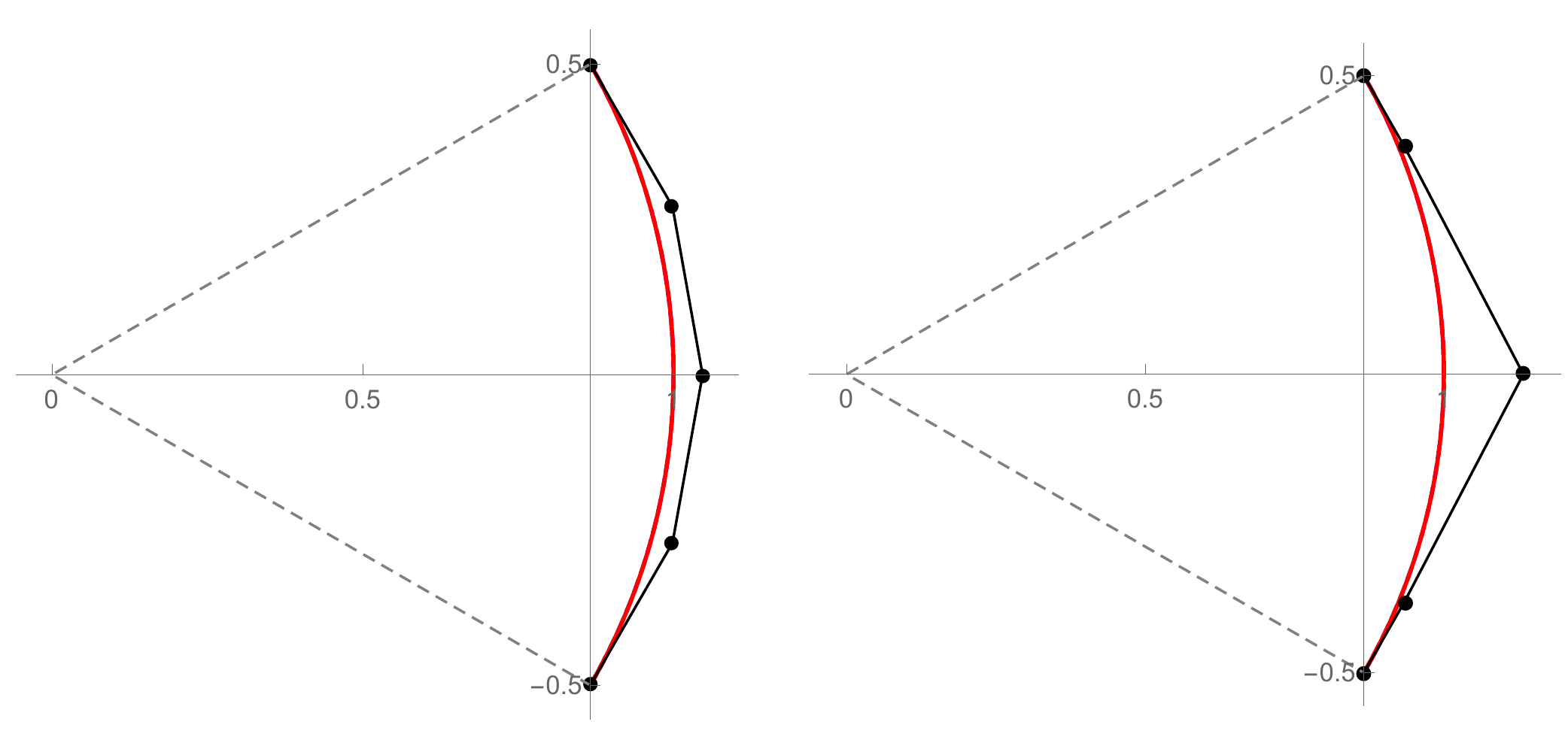}
\caption{Two quartic ${\mathcal G}^{0}$ interpolants of circular arc of inner angle $2\tfrac \pi6$ with the simplified signed error function as in the Figure~\ref{graph_error}. The left one is the best interpolant with the error $2.34778\times 10^{-8}$, the right one has the error $4.01760\times 10^{-5}$.}\label{admissible}
\end{center}
\end{figure}

Numerical computations reveal that the function $f$ has five real zeros for $\varphi\approx 0.9188$, six real zeros for $\varphi< 0.9188$, and only four real zeros for $\varphi> 0.9188$. All real zeros of $f$ induce a quartic ${\mathcal G}^{0}$ interpolant which have an alternating simplified signed error function, i.e., its error function has the same shape as the error function $\chi$ of the best interpolant (only the amplitude can vary). We proved that the largest negative zero of $f$ induces the best interpolant of the unit circular arc with the inner angle $2\varphi$. Examples show that the smallest positive zero of $f$ induces the best interpolant of the unit circular arc with the inner angle $2(\pi-\varphi)$; in this case all three control points $b_1$, $b_2$, $b_3$ of the interpolant lie left from the control points $b_0$ and $b_4$. Most of the remaining zeros induce a non admissible interpolant, i.e., an interpolant with self intersections. The reason why interpolant has self intersections is that some of the control points $b_1$, $b_2$, $b_3$ of an interpolant lie left and some right from the control points $b_0$ and $b_4$.
Numerical computations show that for $\varphi<0.6772$, the second largest negative zero induces an admissible interpolant (all three control points $b_1,b_2,b_3$ are right of the control points $b_0$ and $b_4$) but the corresponding error function has a larger amplitude than the error function of the best interpolant (see Figure \ref{admissible}).

\section{Conclusion}\label{conc}

In this paper we presented an interpolation of a circular arc given by an inner angle $2\varphi$ not greater then $\pi$, where both boundary points of the arc are interpolated. Our method works well in the parabolic case, where for every $\varphi$ we get only one candidate for the best interpolant, and also in the cubic case, where for every $\varphi$ we get only one admissible candidate. In the quartic case we get more candidates and the analysis to figure out which candidate is the best one is quite demanding. Our method could be applied for interpolation of a circular arc by higher order polynomials, but it seems that it is very hard to prove which candidate is the best one. Maybe the method can be used for some particular cases, like half circular arc or quarter circular arc.

{\noindent \sl Acknowledgments.}{ The author would like to thank Emil \v{Z}agar for many useful discussions.

Research on this paper was supported in part by the program P1-0292 and the grants J1-8131 and J1-7025 from ARRS, Republic of Slovenia.}

\end{document}